\newcommand\webcite[1]{\texttt{\def~{\~{}}#1}}

\documentclass{article}

\usepackage{amsmath,amsfonts,amsthm,epsfig,color}
\newtheorem{theorem}{Theorem}
\newtheorem{lemma}[theorem]{Lemma}

\newtheorem{corollary}[theorem]{Corollary}

\def\TT{{\mathbb T}}
\def\LP{{LV}}

\def\epsilon{{\varepsilon}}
\renewcommand\star{{*}}

\newcommand\E{{\mathop{\mathbb E{}}\nolimits}}
\renewcommand\Pr{{\mathop{\mathbb P{}}\nolimits}}

\def\Prone{{\widetilde\Pr}}
\def\cA{{\mathcal A}}
\def\cB{{\mathcal B}}
\def\Pow{{\mathcal P}}
\def\Z{{\mathbb Z}}

\begin{document}
\title{A short proof of the Harris-Kesten Theorem}

\author{B\'ela Bollob\'as\thanks{Department of Mathematical Sciences,
University of Memphis, Memphis TN 38152, USA}
\thanks{Trinity College, Cambridge CB2 1TQ, UK}
\thanks{Research supported in part by NSF grant ITR 0225610 and DARPA grant
F33615-01-C-1900}
\thanks{Research partially undertaken during a visit to the Forschungsinstitut f\"ur Mathematik,
ETH Z\"urich}
\and Oliver Riordan$^{\dag\S}$%
\thanks{Royal Society Research Fellow, Department of Pure Mathematics
and Mathematical Statistics, University of Cambridge, UK}}
\maketitle

\begin{abstract}
We give a short proof of the fundamental result that the
critical probability for bond percolation in the planar square lattice
$\Z^2$ is equal to $1/2$. The lower bound was proved by Harris,
who showed in 1960 that percolation does not occur
at $p=1/2$. The other, more difficult, bound
was proved by Kesten, who showed in 1980 that percolation does occur for any $p>1/2$.
\end{abstract}

\section{Introduction}
Let us set the scene by recalling some basic notions of percolation theory, in
the very special context we shall study here;
for general background, we refer the reader to
Kesten~\cite{Kbook}, Chayes and Chayes~\cite{CC} and Grimmett~\cite{Grimmett}.
Let $\Z^2$ be the planar square lattice, i.e., the graph whose vertices are the points of $\Z^2$,
in which vertices at Euclidean distance $1$ are joined by an edge.
A {\em bond percolation measure} on $\Z^2$, or any other graph, is a probability
measure on the space of assignments of a {\em state}, namely {\em open} or {\em closed},
to each edge $e\in E(\Z^2)$ (with the usual
$\sigma$-field of measurable events). Most of the time we shall consider the product
measure $\Pr_p=\Pr_p^{\Z^2}$, in which the states of the edges are independent, and every edge is open 
with probability $p$. 

An {\em open cluster} is a maximal connected subgraph of $\Z^2$ all of whose edges
are open. We write $C_v$ for the open
cluster containing a given vertex $v\in \Z^2$. Thus a vertex $w$ lies in $C_v$
if and only if $w$ can be reached from
$v$ by an {\em open path}, i.e., a path in $\Z^2$ all of whose edges are open.

The fundamental question of percolation theory is `when does
percolation occur', i.e., `for which $p$ is there
an infinite open cluster'?
Of course, this question can be, and is, asked in a wide variety of
contexts. Here we shall consider only the particular case of {\em bond percolation in $\Z^2$} described
above. Note that the question makes sense: if $E_\infty$ is the event that there is an infinite open
cluster, then by Kolmogorov's $0/1$-law, $\Pr_p(E_\infty)$ is $0$ or $1$ for any $p$.

Writing $|C_v|$ for the number of vertices of $C_v$, let
\[
 \theta(p)=\Pr_p(|C_0|=\infty),
\]
where $0=(0,0)$ is the origin.
Now $|C_0|=\infty$ implies $E_\infty$, while $E_\infty$ is the countable union of the events
$|C_v|=\infty$, $v\in \Z^2$. From translational
invariance we have $\Pr_p(|C_v|=\infty)=\theta(p)$ for all $v$,
so $\Pr_p(E_\infty)$ is $0$ if $\theta(p)=0$ and $1$ if $\theta(p)>0$;
thus, the question of when percolation occurs is precisely the question `when is $\theta(p) > 0$?'
It is easy to see
that $\theta(p)$ is an increasing function of $p$, so there is a `critical probability'
\[
 p_H = p_H(\Z^2) = \inf\{p: \theta(p)>0\} = \inf\{p: \Pr_p(E_\infty)>0\}
 = \inf\{p: \Pr_p(E_\infty)=1\}.
\]
Here, following Welsh (see~\cite{SW}), the $H$ is in honour of Hammersley,
although it is more common to write $p_c$ for $p_H$.
In simple cases such as the present one, it is easy to see that $0<p_H<1$.
The problem of investigating $p_H$ in a variety of contexts
was posed by Broadbent and Hammersley~\cite{BH} in 1957. Hammersley~\cite{H2,H4,H5}
proved general results implying in particular that $0.35<p_H(\Z^2)<0.65$.
The first major
progress was due to Harris~\cite{Harris}, who proved in 1960 that $p_H(\Z^2)\ge 1/2$.
His proof makes use of the `self-duality' of $\Z^2$, and is highly non-trivial.

For many years it was believed that $p_H=1/2$ for bond percolation in $\Z^2$;
see, for example, Sykes and Essam~\cite{SykesEssamPRL}.
However, it was only in 1980, twenty years after Harris proved that $p_H\ge 1/2$,
that Kesten~\cite{Kesten1/2} proved this conjecture, following significant progress
by Russo~\cite{Russo} and Seymour and Welsh~\cite{SW}.

The proofs of Harris and Kesten
are beautiful, but rather complicated and long. Other proofs of
their results have since been developed, for example that due
to Zhang (see~\cite{Grimmett}), and very much more general
results are now known; however, none of the proofs is very simple.
Here we shall give short proofs of the Harris-Kesten results
for $p_H(\Z^2)$, the critical probability for bond percolation
in $\Z^2$.
Our methods are applicable in other contexts as well;
indeed, they were developed in~\cite{Voronoi} to prove that the critical
probability for random Voronoi percolation in the plane is $1/2$.

It has been known for a long time that
Kesten's result follows easily once one can show that a certain
`rectangle crossing' event undergoes a `sharp transition', in
that for any $p>1/2$ its probability is close to $1$ whenever the
rectangle is large enough.
The key to our method is a simple deduction of this result
from a general result of probabilistic combinatorics,
due to Friedgut and Kalai~\cite{FK}.
There are two other components to the proof. For the first,
we use the strategy of Russo~\cite{Russo} and Seymour and Welsh~\cite{SW},
but give a shorter proof.
For the second, deducing percolation from a large rectangle-crossing
probability, we give three arguments, one that `happens to come out'
because of how the numbers work out, a well-known, elegant `renormalization' argument,
and a more recent argument that is likely to be useful in other contexts.

The rest of the paper is organized as follows. In the next section we present
two results from probabilistic combinatorics; with the exception of these
results, the arguments to come will be self-contained.
In Section~\ref{sec_Dual} we discuss the `self-duality' of $\Z^2$.
Harris' Theorem will be proved in Section~\ref{sec_H}, the key application
of the Friedgut-Kalai result will be given in Section~\ref{sec_S}, and three
ways of deducing Kesten's Theorem in Section~\ref{sec_K}.
In the last section
we briefly discuss possible extensions.

\section{Preliminaries}\label{sec_prelim}

The first result we shall need
is a very simple but fundamental lemma due to Harris~\cite{Harris}.
Let $X$ be a fixed ground set with $N$ elements, and let $X_p$ be a random subset
of $X$ obtained by selecting each $x\in X$ independently with probability $p$.
For a family $\cA\subset \Pow(X)$ of subsets of $X$, let $\Pr_p^X(\cA)$ be the probability that $X_p\in\cA$.
In this context, $\cA$ is {\em increasing} if $A\in \cA$ and $A\subset B\subset X$
imply $B\in \cA$. For example, $X$ might be a finite set of edges of $\Z^2$, and $X_p$ the subset
of $X$ consisting of the open edges. Then an event is increasing if
it is preserved when the states of one or more edges are changed
from closed to open.
 
\begin{lemma}\label{l_harris}
If $\cA$, $\cB\subset \Pow(X)$ are increasing, then for any $p$ we have
\[
 \Pr_p^X(\cA\cap \cB)\ge \Pr_p^X(\cA) \Pr_p^X(\cB).
\]
\end{lemma}

In other words, increasing events are positively correlated.
Once one thinks of the statement, the proof turns out to be very simple, using induction on $N$.
Taking complements, Lemma~\ref{l_harris}
also applies to two {\em decreasing} events, defined in the obvious way;
similarly, an increasing event and a decreasing one are negatively correlated.
The extension to infinite product spaces is immediate; we shall not need it here.
Harris discovered this lemma precisely in the context of percolation,
although it is perhaps one of the most basic results of combinatorial
probability. In the latter context Harris' Lemma was rediscovered by
Kleitman~\cite{Kleitman};
a series of important generalizations then followed, culminating in the very
general `Four-functions Theorem' of Ahlswede and Daykin~\cite{AD-FFT}.

The second result we shall need is a sharp-threshold theorem of Friedgut and Kalai~\cite{FK},
which is itself a simple consequence of a result of Kahn, Kalai and Linial~\cite{KKL} (see also~\cite{BKKKL})
concerning
the influences of coordinates in a product space.

With the same set-up as above, $\cA\subset \Pow(X)$ is {\em symmetric} if there is a permutation group acting transitively
on $X$ whose action on $\Pow(X)$ preserves $\cA$.
In our notation, the result of Friedgut and Kalai we shall need is as follows.

\begin{theorem}\label{th_sharp}
There is an absolute constant $c_1$ such that if
$|X|=N$, 
$\cA\subset \Pow(X)$ is symmetric and increasing,
$0<\epsilon<1/2$, and
$\Pr_p^X(\cA)>\epsilon$, then $\Pr_q^X(\cA)>1-\epsilon$ whenever
\begin{equation}\label{q-p}
 q-p \ge c_1 \frac{\log(1/(2\epsilon))}{\log N}.
\end{equation}
\end{theorem}

\vskip2pt
Having collected the basic results above for future reference, we are ready to start our simple
proof that $p_H(\Z^2)=1/2$.

\section{Self-duality}\label{sec_Dual}

Throughout this paper we shall work with `open crossings of rectangles', a 
standard concept of planar percolation.
We shall identify a rectangle $R=[a,b]\times [c,d]$, where $a<b$ and $c<d$ are integers,
with an induced subgraph of $\Z^2$. This subgraph includes all vertices and edges in the interior
and boundary of $R$.
If $b-a=k$ and $d-c=\ell$ then we call $R$ a {\em $k$ by $\ell$ rectangle}.

We write $H(R)$ for the event that there is a {\em horizontal open crossing} of $R$,
i.e., a path from the left side of $R$ to the right side consisting entirely of open edges
of $R$.
Similarly, we write $V(R)$ for the event that there is a vertical open crossing
of $R$.
Note that for disjoint rectangles $R$ and $R'$ (a rectangle includes its boundary)
each of the events $H(R)$, $V(R)$ is independent of each of the events $H(R')$, $V(R')$.

The {\em dual} $L^\star$ of the lattice $L=\Z^2$ is the usual planar dual of $\Z^2$ considered
as a graph drawn in the plane: there is a vertex $v\in L^\star$ corresponding to each face of $\Z^2$,
i.e., to each square $[x,x+1]\times [y,y+1]$, $x,y\in \Z$. 
It is customary to take $v=(x+1/2,y+1/2)$.
There is one edge $e^\star$ of $L^\star$ corresponding to each edge $e$ of $L$;
this edge joins the two vertices of $L^\star$
corresponding to the faces of $L$ in whose boundary $e$ lies. Of course, $L^\star$ may be defined analogously
for any plane graph. When $L=\Z^2$, then $L^\star=(1/2,1/2)+\Z^2$ is isomorphic
to $L$, as shown below.
\begin{figure}[htb]
 \[\epsfig{file=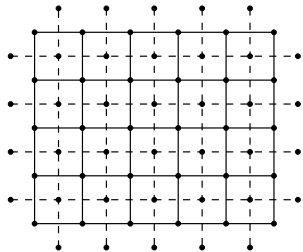,height=1.3in}\]
\caption{Portions of the lattice $L=\Z^2$ (solid lines) and the isomorphic dual lattice $L^\star$ (dashed lines).}
\label{fig_selfdual}
\end{figure}

The {\em horizontal dual}, or simply the {\em dual}, of a rectangle $R=[a,b]\times[c,d]$
is the rectangle
$R^h=[a+1/2,b-1/2]\times [c-1/2,d+1/2]$ in $L^\star$.
Defining an edge $e^\star$ of $L^\star$ to be open if and only if $e$ is closed,
let $V^\star(R^h)$ be the event that there is a vertical crossing of $R^h$ by open edges of $L^\star$.

The next lemma is key to the analysis of bond percolation in $\Z^2$. Often,
it is stated as `obvious' and no formal proof is given. While the result is indeed
obvious, it is not entirely
trivial to prove. However, a short proof is possible, and, as it
happens, needs no topology. For a topological proof of a related result
see Kesten~\cite[pp. 386--392]{Kbook}.

\begin{lemma}\label{l_oneway}
Let $R$ be a rectangle in $L=\Z^2$. Whatever the states of the edges in $R$,
exactly one of the events $H(R)$ and $V^\star(R^h)$ holds.
\end{lemma}

\begin{proof}
To avoid fractions, for this proof let $L=(0,2)+4{\mathbb Z}^2$ and
$L^*=(2,0)+4{\mathbb Z}^2$; also, let $R\subset L$ have vertex set
$\{(4i, 4j+2): \, 0\le i \le k,\ 0\le j \le \ell-1 \}$, so that the
vertex set of $R^h$ is $\{(4i+2, 4j): \, 0\le i \le k-1, \ 0\le j \le
\ell \}$. The open edges of $R$ form a graph $G$, and those of
$R^h$ form a graph $G^h$. Our task is to show that either $G$ contains a
left-right path,
i.e., a path from a vertex $(0, 4j+2)$ to a vertex
$(4k, 4h+2)$, or else $G^h$ contains a top-bottom path, but not both.
In proving this assertion, we may assume that $G$ contains all
$2(\ell -1)$ edges on the left and right sides,
and $G^h$ contains all
$2(k-1)$ edges at the top and bottom, as in Figure~\ref{fig_oneway}. Note that
$G$ and $G^h$ are plane graphs with every edge a straight-line
segment of length 4.

\begin{figure}[htb]
 \centering
 \input{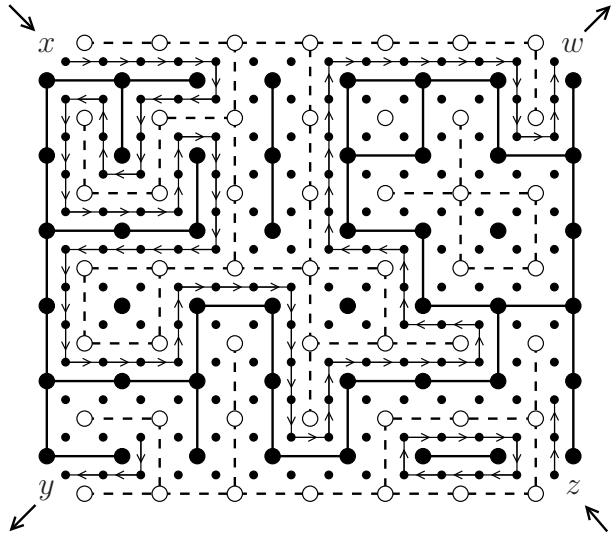}
 \caption{The graphs $G$ (large solid circles and solid lines) and $G^h$ (empty circles
and dashed lines), and parts of the oriented graph $M$ (small circles and arrows). The part
of $M$ shown includes all of the path $P$ and a cyclic component.}\label{fig_oneway}
\end{figure}

We shall construct a third graph, $M$, which is in the `middle',
between $G$ and $G^h$.
Let $R_{\rm odd}$ be the $2k \times
2\ell$ rectangle with vertex set $\{ (2i+1, 2j+1):
\, 0\le i \le 2k-1, \, 0\le j \le 2\ell-1 \, \}$, and let $M$ be the
subgraph of $R_{\rm odd}$ formed by the edges not meeting either
an edge of $G$ or an edge of $G^h$. Orient every edge of $M$
so that $G$ is on the right and $G^h$ is on the left.
Then every vertex of $M$ has one edge coming in and one going out,
apart from exactly four vertices of
degree one, $x, y, z$ and $w$, say, as in Figure~\ref{fig_oneway}. Thus
the component of $x$ in $M$ is a path $P$ ending at one of $y$ and $w$.
If $P$ ends at $w$ then the edges of $G$
to the right of $P$ form a connected subgraph of $G$ meeting the
left side {\em and} the right side (as in Figure~\ref{fig_oneway}), and if $P$
ends at $y$ then the edges of $G^h$ to the left of $P$ form a
connected subgraph of $G^h$ meeting the top {\em and} the bottom.
Since a connected subgraph meeting two opposite sides contains a
path meeting the same sides, we have shown that at least
one of $H(R)$ and $V^\star(R^h)$ holds.

That $H(R)$ and $V^\star(R^h)$ cannot both hold is immediate from the Jordan
Curve
Theorem; however, only the easy part of this theorem is needed.
Indeed, if the path $P$ ends at $w$ then $P$ may be completed to a 
simple closed curve $\overline P$ by coming back from $w$ to $x$ over the
top of $R^h$. This curve winds once around any top vertex of $R^h$,
and zero times around any bottom vertex. As the winding number of
$\overline P$ is constant off $\overline P$, any top-bottom path in $G^h$ must
cross $\overline P$, and hence cross an edge of $M$, contradicting
the definition of $M$. Similarly, if $P$ ends at $y$ there
is no left-right path in $G$, completing the proof.
\end{proof}

As the edges of $L^\star$ are open independently with probability $1-p$,
Lemma~\ref{l_oneway} has the following immediate corollary concerning
rectangles in the original lattice $L=\Z^2$.

\begin{corollary}\label{c_sumto1}
If $R$ is a $k$ by $\ell-1$ rectangle in $\Z^2$ and $R'$ is a $k-1$ by $\ell$ rectangle,
then
\[
 \Pr_p(H(R))+\Pr_{1-p}(V(R')) =1.
\]
\end{corollary}

If $k=\ell=n+1$, then $R'$ is $R$ rotated by $90$ degrees.
Thus $\Pr_p(H(R))=\Pr_p(V(R'))$, so Corollary~\ref{c_sumto1}
implies the following essential consequence of the self-duality of $\Z^2$.

\begin{corollary}\label{c_half}
If $R$ is an $n+1$ by $n$ rectangle then $\Pr_{1/2}(H(R))=1/2$.
Thus, if $S$ is an $n$ by $n$ square, then
\[
 \Pr_{1/2}(V(S)) = \Pr_{1/2}(H(S)) \ge 1/2.
\]
\end{corollary}

It is easy to think that Corollary~\ref{c_half} shows that $p_H=1/2$.
Although self-duality is of course the reason `why'
$p_H=1/2$, a rigorous deduction is far from easy, and took twenty years to accomplish.

\section{A short proof of Harris' Theorem}\label{sec_H}

In this section we present a short proof of Harris' 1960 result that $\theta(1/2)=0$,
and hence $p_H=p_H(\Z^2)\ge 1/2$. The techniques are considerably
simpler than Harris' original ones. The strategy of the proof is that
of Russo~\cite{Russo} and Seymour and Welsh~\cite{SW}; however, our proof
of the key intermediate result is considerably
shorter than theirs.

\begin{lemma}\label{l_X}
Let $R=[0,m]\times [0,2n]$, $m\ge n$, be an $m$ by $2n$ rectangle. Let $X(R)$ be the event that
there are paths $P_1$ and $P_2$ of open edges, where $P_1$ crosses $S=[0,n]\times [0,n]$ from top to bottom,
and $P_2$ lies inside $R$ and joins some vertex on $P_1$ to some vertex on the right-hand side of $R$.
Then $\Pr_p(X(R)) \ge \Pr_p(H(R))\Pr_p(V(S))/2$.
\end{lemma}

\begin{figure}[htb]
 \centering
 \input{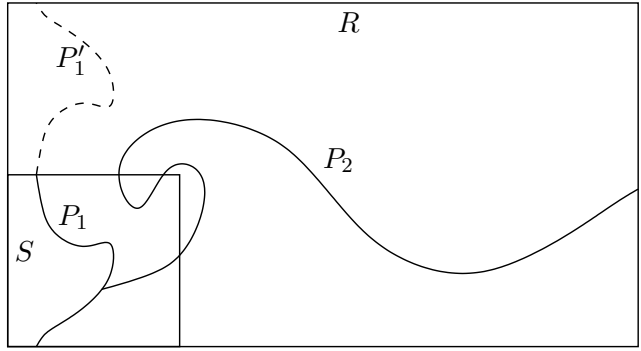}
\caption{A rectangle $R$ and square $S$ inside it, drawn with paths (solid curves) whose
presence as open paths would imply $X(R)$.}\label{fig_XRdef}
\end{figure}

\begin{proof}
Suppose that $V(S)$ holds, so there is a path $P$ of open edges crossing $S$ from top to bottom. Note
that any such $P$ separates $S$ into two pieces, one to the left of $P$ and one to the right.
The proof of Lemma~\ref{l_oneway} shows that when $V(S)$ holds,
one can define the left-most such path $P$, $\LP(S)$, in such a way that the event
$\{\LP(S) = P_1\}$ does not depend on edges to the right of $P_1$.

It is easy to see that for any possible value $P_1$ of $\LP(S)$ we have
\[
 \Pr_p(X(R)\mid \LP(S)=P_1)\ge \Pr_p(H(R))/2:
\]
let $P$ be the (not necessarily open)
path formed by the union of $P_1$ and its reflection $P_1'$ in the line $y=n$; see Figure
\ref{fig_XRdef}. This path crosses
$[0,n]\times [0,2n]$ from top to bottom. With (unconditional) probability $\Pr_p(H(R))$ there is a path
$P_3$ of open edges crossing $R$ from right to left - this path must meet $P$. By symmetry, the (unconditional)
probability that some such path first meets $P$ at a point of $P_1$ is at least $\Pr_p(H(R))/2$.
Hence, the event $Y(P_1)$ that there is a path $P_2$ in $R$ to the right of $P$ joining some point of $P_1$ to the
right-hand side of $R$ has probability at least $\Pr_p(H(R))/2$.
But $Y(P_1)$ depends only on edges to the right of $P$.
All such edges in $S$ are to the right of $P_1$ in $S$. As the states of these edges are independent
of $\{\LP(S) = P_1\}$, we have
\[
 \Pr_p(Y(P_1)\mid \LP(S)=P_1)=\Pr_p(Y(P_1))\ge \Pr_p(H(R))/2.
\]
But $\LP(S)=P_1$ and $Y(P_1)$ imply
$X(R)$, so we have $\Pr_p(X(R)\mid \LP(S)=P_1)\ge \Pr_p(H(R))/2$. Finally, the event $V(S)$ is a disjoint union
of events of the form $\{\LP(S)=P_1\}$, so $\Pr_p(X(R)\mid V(S))\ge \Pr_p(H(R))/2$.
\end{proof}

The following immediate corollary of Lemma~\ref{l_X} is standard, but is usually proved in a different way
(by the methods of Russo~\cite{Russo} or Seymour and Welsh~\cite{SW}; see~\cite{Grimmett}).

\begin{corollary}\label{c_long}
Let $\rho>1$ be a fixed integer. There is a constant $c_2(\rho)>0$ depending only on $\rho$ such that
for any $2\rho n$ by $2n$ rectangle $R$ we have $\Pr_{1/2}(H(R))\ge c_2(\rho)$.
\end{corollary}

Note that we take $\rho$
to be an integer only because we want all our rectangles to have integer
coordinates.

\begin{proof}
Writing $h_{m,n}$ for $\Pr_{1/2}(H(R))$, where $R$ is an $m$ by $n$ rectangle, we claim that
for $m\ge n$ we have $h_{2m-n,2n}\ge h_{m,2n}^2/2^5$. Applying this repeatedly, starting from
$h_{2n,2n}\ge 1/2$ (by Corollary~\ref{c_half}),
we obtain the result.

To see the claim, consider the rectangles $R=[0,m]\times [0,2n]$, $R'=[n-m,n]\times [0,2n]$
and the square $S=[0,n]\times [0,n]$ in their intersection, as shown in Figure~\ref{fig_RSR}.

\begin{figure}[htb]
\centering
\input{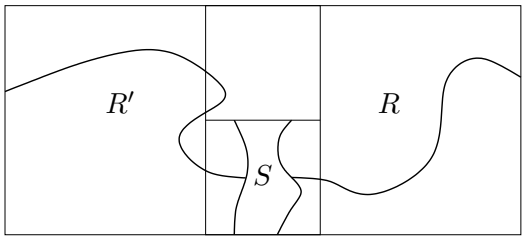}
\caption{The overlapping rectangles $R$ and $R'$ with the square $S$ in their intersection. The paths drawn show
that $X(R)$ holds, as well as the reflected equivalent for $R'$. If $H(S)$ also holds, then so does $H(R\cup R')$.}
\label{fig_RSR}
\end{figure}

If $E_1=X(R)$, the corresponding
(i.e., horizontally reflected) event $E_2$ for $R'$, and $E_3=H(S)$ all hold, then so does $H(R\cup R')$,
using only the fact that any horizontal crossing of $S$ meets any vertical crossing. (See Figure~\ref{fig_RSR}.)
But the $E_i$ are increasing events, so, by Lemma~\ref{l_harris},
\[
 \Pr_{1/2}(H(R\cup R')) \ge \Pr_{1/2}(E_1)\Pr_{1/2}(E_2)\Pr_{1/2}(E_3) = \Pr_{1/2}(X(R))^2\Pr_{1/2}(H(S)).
\]
By Corollary~\ref{c_half}, $\Pr_{1/2}(V(S))=\Pr_{1/2}(H(S))\ge 1/2$.
Applying Lemma~\ref{l_X}, it follows that
\begin{eqnarray*}
 \Pr_{1/2}(H(R\cup R')) &\ge& \Pr_{1/2}(H(R))^2\Pr_{1/2}(V(S))^2\Pr_{1/2}(H(S))/4 \\
 &\ge& \Pr_{1/2}(H(R))^2/2^5.
\end{eqnarray*}
This proves the claim and thus the corollary.
\end{proof}

As noted by Russo~\cite{Russo} and Seymour and Welsh~\cite{SW},
Harris' result is an easy consequence of Corollary~\ref{c_long}.

\begin{theorem}\label{th_Harris}
For bond percolation in $\Z^2$, $\theta(1/2)=0$.
\end{theorem}

\begin{proof}
Let $c_2=c_2(3)$ be the absolute constant given by Corollary~\ref{c_long}.
By Corollary~\ref{c_long} and Lemma~\ref{l_harris}, if we arrange
four $6n$ by $2n$ rectangles to form a `square annulus' as in Figure~\ref{fig_ann},
\begin{figure}[htb]
 \[\epsfig{file=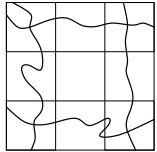,height=0.8in}\]
\caption{Four rectangles forming a square annulus.}
\label{fig_ann}
\end{figure}
then with probability at least $c_2^{\,4}$
each rectangle is crossed the long way by an open path. The union of any such paths contains a cycle
surrounding the centre of the annulus; see Figure~\ref{fig_ann}.
For $k\ge 1$, let $A_k$ be the square annulus centred on the origin with
inner and outer radii $4^k$ and $3\times 4^k$. Then the $A_k$ (including their boundaries) are disjoint.
Using independence of disjoint regions, at $p=1/2$
the probability that the point $(1/2,1/2)$ is {\em not} surrounded by a cycle of open edges
is at most
\[
 \prod_{k=1}^\infty \Pr_{1/2}(\hbox{$A_k$ contains no open cycle}) \le \prod_{k=1}^\infty (1-c_2^{\,4}) =0.
\]
Hence, at $p=1/2$, the probability that $(1/2,1/2)$ is in an infinite
open cluster in the dual lattice is zero.
Equivalently, at $p=1/2$ the probability that there is an infinite open cluster is $0$.
\end{proof}

\section{Crossing rectangles with high probability}\label{sec_S}

Surprisingly, although it was twenty years after Harris proved that $p_H(\Z^2)\ge 1/2$ that
Kesten found an ingenious proof that $p_H(\Z^2)=1/2$, one can use the
results mentioned in Section~\ref{sec_prelim}
to deduce Kesten's Theorem from Corollary~\ref{c_long}
with very little additional work. The main idea is to
find a way of using Theorem~\ref{th_sharp} to show that large rectangles
can be crossed with probability close to $1$.

\begin{lemma}\label{l_long}
Let $p>1/2$ and an integer $\rho>1$ be fixed.
There are constants $\delta=\delta(p)>0$ and $n_0=n_0(p,\rho)$ such that,
if $n\ge n_0$ and
$R_n$ is a $4\rho n$ by $4n$ rectangle in $\Z^2$, then
$\Pr_p(H(R_n))\ge 1-n^{-\delta}$.
\end{lemma}

\begin{proof}
To simplify the notation we shall prove the lemma in the case $\rho=3$, which
immediately implies the case $\rho=2$; these are the only cases we shall use.
The full result can be proved in the same way; it also follows from the case $\rho=2$, say,
by combining crossings as in Figure~\ref{fig_4by1}.

Fix $p>1/2$. From Corollary~\ref{c_long}, applied with $\rho=7$, 
there is an absolute constant $c_2$, $0<c_2<1/2$, such that
$\Pr_{1/2}(H(R))\ge c_2$ for any $14n$ by $2n$ rectangle $R$.
If $H(R)$ were a symmetric event then we could just apply Theorem~\ref{th_sharp},
but $H(R)$ is not symmetric, so we shall introduce symmetric auxiliary events.

For $n\ge 3$, let $\TT_n$ be the $n$ by $n$ discrete torus, i.e., the graph $C_n\times C_n$, which may be obtained
from $\Z^2$ by identifying all pairs of vertices for which the corresponding coordinates
are congruent modulo $n$. Thus $\TT_n$ has $n^2$ vertices and $2n^2$ edges.
For $1\le k,\ell\le n-2$, a {\em $k$ by $\ell$ rectangle} $R$ in $\TT_n$ is an induced subgraph of $\TT_n$
corresponding to a $k$ by $\ell$ rectangle $R'=[a,a+k]\times [b,b+\ell]$ in $\Z^2$.
Note that our rectangles in the torus are always too small to `wrap around', so the induced subgraph of $\TT_n$
is isomorphic to the corresponding induced subgraph of $\Z^2$.

We shall take the edges of $\TT_n$ to be open independently with probability $p$, writing
$\Pr_p^{\TT_n}$ for the corresponding probability measure. Most of the time, we shall suppress the dependence
 on $n$.

Let $E_n$ be the event that $\TT_{16n}$ contains {\em some} $14n$ by $2n$
rectangle with a horizontal crossing by open edges, or some $2n$ by
$14n$ rectangle with a vertical crossing. Note that $E_n$ is
symmetric as a subset of $\Pow(X)$, where $X$ is the edge set of
$\TT_{16n}$, which has cardinality $N=512n^2$.

Considering one fixed $14n$ by $2n$ rectangle $R$ in the torus, which we may
identify with a corresponding rectangle in $\Z^2$, and writing $\TT$ for $\TT_{16n}$, we have
\[
 \Pr_{1/2}^{\TT}(E_n)\ge \Pr_{1/2}^\TT(H(R)) = \Pr_{1/2}(H(R)) \ge c_2.
\]
Let $\delta=(p-1/2)/(64c_1)$, where $c_1$ is the constant in Theorem~\ref{th_sharp},
and set $\epsilon=n^{-128\delta}$. As $\delta$ depends only on $p$,
there is an $n_0=n_0(p)$ such that $\epsilon < c_2\le 1/2$
for all $n\ge n_0$.
Now
\[
 p-\frac{1}{2} = 64c_1\delta = c_1\frac{\log(1/\epsilon)}{\log(n^2)}
 > c_1\frac{\log(1/(2\epsilon))}{\log(512n^2)}.
\]
Hence, by Theorem~\ref{th_sharp}, as $\Pr_{1/2}^\TT(E_n)\ge c_2 > \epsilon$, we have
\begin{equation}\label{e2}
 \Pr_{p}^\TT(E_n) \ge 1-\epsilon = 1-n^{-128\delta}
\end{equation}
for all $n\ge n_0$.

Clearly, there are $12n$ by $4n$ rectangles $R_1,\ldots,R_{64}$ covering $\TT=\TT_{16n}$ 
such that
any $14n$ by $2n$ rectangle $R$ in $\TT$ crosses one of the $R_i$, in the sense
that the intersection of $R$ and $R_i$ is a $12n$ by $2n$ subrectangle of $R_i$. For example,
we may take the $12n$ by $4n$ rectangles $R_i$ whose bottom-left coordinates are all possible multiples of $2n$.
Similarly, there are $4n$ by $12n$ rectangles $R_{65},\ldots,R_{128}$ 
so that any $2n$ by $14n$ rectangle in $\TT$ crosses one of these from top to bottom.

It follows that if $E_n$ holds, then so does one of the events $E_{n,i}$, $i=1,\ldots,128$,
where each $E_{n,i}$ is the event that $R_i$ is crossed
the long way by an open path. Thus $E_n^c$, the complement of $E_n$, contains the intersection
of the $E_{n,i}^c$.

By Lemma~\ref{l_harris}, applied to the product measure $\Pr_p^\TT$, for each $i$
the decreasing events $E_{n,i}^c$ and $\cap_{j<i} E_{n,j}^c$ are positively correlated.
Hence
\[
 \Pr_p^\TT(E_n^c) \ge \Pr_p^\TT\left(\bigcap_{i=1}^{128} E_{n,i}^c\right)
 \ge \prod_{i=1}^{128}\Pr_p^\TT(E_{n,i}^c) = \Pr_p^\TT(E_{n,1}^c)^{128}.
\]
Thus, from \eqref{e2}, for $n\ge n_0$ we have
\[
 \Pr_p^\TT(E_{n,1}^c) \le \Pr_p^\TT(E_n^c)^{1/128} \le n^{-\delta},
\]
so $\Pr_p^\TT(E_{n,1})\ge 1-n^{-\delta}$.

Now $E_{n,1}$ is the event that there is a horizontal open
crossing of a fixed $12n$ by $4n$ rectangle $R$ in the torus, which we may
identify with a corresponding rectangle in $\Z^2$.
Thus
\begin{equation}\label{longhigh}
 \Pr_p(H(R)) = \Pr_p^{\TT_{16n}}(H(R)) \ge 1-n^{-\delta}
\end{equation}
whenever $R$ is a $12n$ by $4n$ rectangle in $\Z^2$ and $n$ is large enough,
completing the proof of the lemma.
\end{proof}

\section{Three ways of deducing Kesten's Theorem}\label{sec_K}

Our first deduction makes use of the quantitative form of Lemma~\ref{l_long}.
Recall that $E_\infty$ is the event that there is an infinite open cluster,
and that $\Pr_p(E_\infty)>0$ implies that $\Pr_p(E_\infty)=1$ and $\theta(p)>0$.

\begin{theorem}\label{th_K}
For bond percolation in $\Z^2$, if $p>1/2$ then $\Pr_p(E_\infty)=1$.
\end{theorem}

\begin{proof}[Proof of Theorem~\ref{th_K} -- first version]
Fix $p>1/2$.
Let $\delta=\delta(p)$ and $n_0=n_0(p,2)$ be as in Lemma~\ref{l_long}.
Let $m\ge n_0$ be an integer to be chosen below, and set $n=4m$.
For $k=0,1,2,\ldots$, let $R_k$ be a rectangle with bottom-left corner the origin
and side-lengths $2^kn$ and $2^{k+1}n$, where the longer side is vertical if
$k$ is even and horizontal if $k$ is odd; see Figure~\ref{fig_Rk}.

\begin{figure}[htb]
\centering
\input{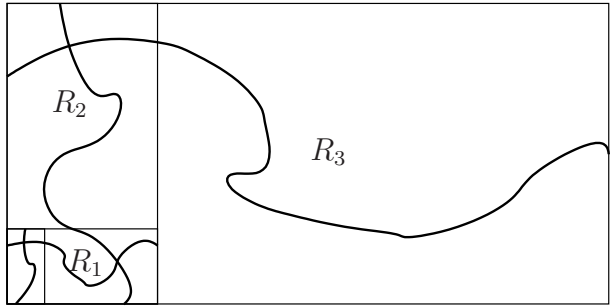}
\caption{The rectangles $R_0$ to $R_3$ ($R_0$ not labelled) drawn with open
paths corresponding to the events $E_k$.}
\label{fig_Rk}
\end{figure}

Let $E_k$ be the event that $R_k$ is crossed the long way by an open path.
Note that any two such crossings of $R_k$ and $R_{k+1}$ must meet, so
if all the $E_k$ hold, then so does $E_\infty$.
If $n$ is large enough then,
by Lemma~\ref{l_long},
\[
 \sum_{k\ge0} \Pr_p(E_k^c)
 \le \sum_{k\ge 0} (2^kn)^{-\delta} = \frac{n^{-\delta}}{1-2^{-\delta}} <1,
\]
so $\Pr_p(E_\infty)\ge \Pr_p(\bigcap_{k\ge0} E_k)>0$.
\end{proof}

Together, Theorems~\ref{th_Harris} and~\ref{th_K} show that $p_H(\Z^2)=1/2$.

\medskip
The argument above depends on the result of a certain
calculation: it is important that the sum of the `error probabilities'
$n^{-\delta}$ as $n$ runs over powers of $2$ is convergent. It might appear
that the basic strategy of the proof thus depends on the serendipitous form of the
bound on $q-p$ in \eqref{q-p}. In fact, this is not the
case: one only needs
the qualitative result that for each $\epsilon>0$, the bound
tends to zero as $n\to\infty$.
Starting from this weaker result, the argument of Section~\ref{sec_S}
implies the following qualitative form of Lemma~\ref{l_long}.

\begin{lemma}\label{l_long_weak}
Let $p>1/2$ be fixed.
If $R_n$ is a $3n$ by $n$ rectangle in $\Z^2$, then
$\Pr_p(H(R_n))\to 1$ as $n\to\infty$.
\end{lemma}

It is well known that this lemma implies Kesten's Theorem.
We give two arguments.
The first is a `renormalization' argument due to
Aizenman, Chayes, Chayes, Fr\"ohlich and Russo~\cite{ACCFR};
see also Chayes and Chayes~\cite{CC}.

\begin{proof}[Proof of Theorem~\ref{th_K} -- second version]
Let us call a crossing of a rectangle {\em internal}
if it uses no edges in the boundary of the rectangle.
Fixing $p>1/2$ as before,
let $i_n$ be the probability that a $2n$ by $n$ rectangle has an internal
horizontal open crossing. Note that $i_n$ is the probability that
a $2n$ by $n-2$ rectangle has a horizontal open crossing.

Writing $s_n$ for the probability that an $n$ by $n$ square has an internal
vertical open crossing, as an $n$ by $n$ square is
contained in an $n$ by $2n$ rectangle, $s_n\ge i_n$.
Considering three $2n$ by $n$ rectangles overlapping
in two $n$ by $n$ squares as in Figure~\ref{fig_4by1},
by Lemma~\ref{l_harris} the probability that a $4n$ by
$n$ rectangle has an internal horizontal open crossing is at least
$i_n^3s_n^2$. Placing two such rectangles side by side to form a $4n$ by
$2n$ rectangle, the events that each has an internal open crossing
depend on disjoint sets of edges
and are thus independent.
Hence,
\begin{equation}\label{irec}
 i_{2n} \ge 1-\left(1-i_n^3s_n^2\right)^2 \ge 1-\left(1-i_n^5\right)^2.
\end{equation}

\begin{figure}[htb]
\centering
\[\epsfig{file=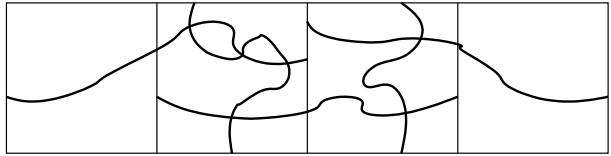}\]
\caption{Three $2n$ by $n$ rectangles and two $n$ by $n$ squares, drawn
with paths guaranteeing a horizontal open crossing of their union}\label{fig_4by1}
\end{figure}

Writing $i_n$ as $1-\epsilon$, from \eqref{irec} we have $i_{2n}\ge 1-25\epsilon^2$,
which is at least $1-\epsilon/2$ if $\epsilon\le 1/50$.
By Lemma~\ref{l_long_weak} there is an $n$ with $i_n\ge 0.98$. (If $n\ge 6$
then a crossing of a fixed $3(n-2)$ by $n-2$ rectangle includes
an internal crossing of a fixed $2n$ by $n$ rectangle.)
It follows that $i_{2^kn}\ge 1-2^{-k}/50$, and $\Pr_p(E_\infty)>0$ follows
as in the first proof of Theorem~\ref{th_K}.
\end{proof}
This method shows that if for
a single value of $n$ we have $i_n> 0.951...$ (a root of
$x=1-(1-x^5)^2$), then percolation occurs.

In fact, arguing as in Chayes and Chayes~\cite{CC}, one can do a little better.
Note that $s_n$ is also the probability that an $n$ by $n$ square has an internal {\em horizontal}
open crossing. Dividing a $2n$ by $n$ rectangle into two squares, if
the rectangle has an internal horizontal open crossing, so do both squares. As the internal
edges of the squares are disjoint, it follows that $i_n\le s_n^2$.
Therefore the first inequality in \eqref{irec} implies
$i_{2n}\ge 1-(1-i_n^4)^2$, and
the value $0.951...$ may be replaced by $0.920...$, a root of
$x=1-(1-x^4)^2$.

\medskip
The second argument requires an even weaker initial bound on $i_n$;
for this
we shall need an observation concerning {\em $k$-dependent percolation}.
A bond percolation measure on $\Z^2$ is {\em $k$-dependent} if
for every pair $S$, $T$ of sets of edges of $\Z^2$ at graph
distance at least $k$, the states (being open or closed) of the edges in $S$ are independent of
the states of the edges in $T$. When $k=1$, the separation
condition is exactly that no edge of $S$ shares a vertex with an edge of $T$.

In static renormalization arguments, $k$-dependent probability measures arise very naturally.
Comparisons between such measures and product measures (or arguments amounting to such comparisons)
have been considered by several authors; see
Liggett, Schonmann and Stacey~\cite{LSS} and the references therein.

\begin{lemma}\label{l_1dep}
There is a $p_0<1$ such that in any $1$-dependent bond percolation measure 
on $\Z^2$ satisfying the
additional condition that each edge is open with probability at least $p_0$, the probability
that $|C_0|=\infty$ is positive.
\end{lemma}

In many contexts the value of $p_0$ is important. The best bound known is
the result that one can take $p_0=0.8639$; this was proved by
Balister, Bollob\'as and Walters~\cite{BBW}, who used it
in the study of random geometric graphs.
For proving Kesten's Theorem,
however, the value of $p_0$ is not important.

Lemma~\ref{l_1dep} is an immediate consequence of the very general main result
of~\cite{LSS} but, in the form above, is trivial from first principles.
Indeed, if $C_0$ is finite,
then it is surrounded by an open cycle
in the dual lattice $L^\star$.
Very crudely, there are at most $\ell3^\ell$ cycles
of length $\ell$ in $L^\star$ surrounding the origin. For each, the corresponding
set $S$ of edges of $L=\Z^2$ contains a subset $S'$ of size at least $|S|/4=\ell/4$
in which any two edges are vertex disjoint.
By $1$-dependence, the states of the edges in $S'$ are independent, so
the probability that all these edges are closed in $L$ (and thus open in $L^\star$)
is at most $(1-p_0)^{\ell/4}$. Combining these observations shows that
$\Pr(|C_0|<\infty) \le \sum_{\ell\ge 4} \ell3^\ell (1-p_0)^{\ell/4}$, which is less than $1$ if $p_0$
is large enough.

\begin{proof}[Proof of Theorem~\ref{th_K} -- third version]
Let $p>1/2$ be fixed, let $p_0<1$ be a constant for which Lemma~\ref{l_1dep} holds,
and set $c=p_0^{1/3}$.
Given a $3n$ by $n$ rectangle $R$, let $S'$ and $S''$ be the two end squares when $R$ is cut into three
squares. Note that $H(R)$ certainly implies $H(S')$ so, by Lemma~\ref{l_long_weak},
\[
 \Pr_p(V(S''))=\Pr_p(V(S'))=\Pr_p(H(S')) \ge \Pr_p(H(R))\ge c
\]
if $n$ is large enough, which we shall assume from now on.

Let $G(R)$ be the event $H(R)\cap V(S')\cap V(S'')$; see Figure~\ref{fig_GR}.

\begin{figure}[htb]
 \[\epsfig{file=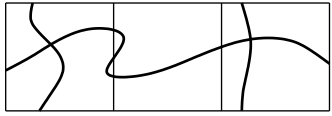,height=0.6in}\]
\caption{A $3n$ by $n$ rectangle $R$ such that $G(R)$ holds.}
\label{fig_GR}
\end{figure}
\noindent

By Lemma~\ref{l_harris},
\[
 \Pr_p(G(R)) \ge \Pr_p(H(R))\Pr_p(V(S'))\Pr_p(V(S'')) \ge c^3=p_0.
\]
Define $G(R')$ similarly for an $n$ by $3n$ rectangle, so $\Pr_p(G(R'))=\Pr_p(G(R))\ge p_0$.

Let us define a $1$-dependent bond percolation measure $\Prone$ on $\Z^2$
as follows: the edge from $(x,y)$ to $(x+1,y)$
is open in $\Prone$ if and only if $G(R)$ holds in $\Pr_p$ for the $3n$ by $n$ rectangle $[2nx,2nx+3n]\times [2ny,2ny+n]$.
Similarly, the edge from $(x,y)$ to $(x,y+1)$ is open in $\Prone$ if and only if $G(R')$ holds in $\Pr_p$ for the $n$
by $3n$ rectangle $[2nx,2nx+n]\times [2ny,2ny+3n]$.
\begin{figure}[htb]
 \[\epsfig{file=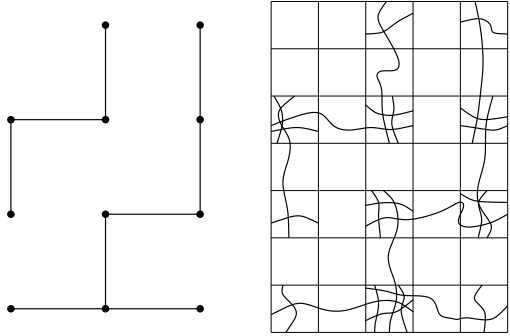,height=1.75in}\]
\caption{A set of open edges in $\Prone$ (left), and corresponding rectangles $R$ drawn with $G(R)$ holding in $\Pr_p$.}
\label{fig_1indep}
\end{figure}

This probability measure is indeed $1$-dependent,
as $G(R)$ depends only on the states of edges in $R$,
and vertex disjoint edges of $\Z^2$ correspond to disjoint rectangles.

By Lemma~\ref{l_1dep}, $\Prone(|C_0|=\infty)>0$.
However, we have defined $G(R)$ in such a way
that a $\Prone$-open path guarantees a corresponding (much longer)
open path in the original bond percolation,
using the fact that horizontal and vertical crossings of a square must meet; see Figure~\ref{fig_1indep}.
Hence, $\Pr_p(E_\infty)\ge \Prone(|C_0|=\infty)>0$,
completing the proof of Theorem~\ref{th_K}.
\end{proof}

\medskip
The argument above works with $2$ by $1$ rectangles, using internal crossings
as in the second proof of Theorem~\ref{th_K}.
Also, it is enough to require a vertical crossing of the left-hand end square of
each rectangle $R$, and a horizontal crossing of the bottom square of each $R'$.
Appealing to $s_n\ge i_n^{1/2}$, as before, to prove percolation it thus suffices to find an $n$
with $i_n^{3/2}\ge p_0$, where $p_0$ is a constant for which Lemma~\ref{l_1dep}
holds. Using the value $p_0=0.8639$ from \cite{BBW}, $i_n\ge 0.907...$
will do.

\section{Extensions}

The arguments above give short proofs of Kesten's Theorem, using Theorem~\ref{th_sharp}
as a key ingredient.
In fact, as we shall show in future work~\cite{ourKesten2},
these arguments, like those of Harris and Kesten,
easily give further results. For example, our method gives
exponential decay of $|C_0|$ when $p<1/2$, showing that $p_T$, the critical probability
at which $\E(|C_0|)$ becomes infinite, is also equal to $1/2$, another result of Kesten~\cite{Kesten1/2}.
Furthermore, although we have
written everything for bond percolation in $\Z^2$, the same method gives similar results in other
contexts; for example, we obtain new, simple proofs
that $p_T=p_H$ for site percolation in the square lattice, and that
both critical probabilities are equal to $1/2$ for site percolation in the triangular lattice.
The basic method used in Section~\ref{sec_S} is much more general, applying to many other two-dimensional
contexts. Indeed, as remarked earlier, it was developed in~\cite{Voronoi} to prove that the critical probability for
random Voronoi percolation in the plane is $1/2$.

\bigskip\noindent
{\bf Acknowledgements.}
Since the first draft of this paper was written we have had stimulating discussions with
Christian Borgs, Jennifer Chayes and Geoffrey Grimmett.
We are grateful to all of them.

\end{document}